\documentclass[11pt]{amsart}%
\usepackage{amsfonts}
\usepackage{amsmath}
\usepackage{amssymb}
\usepackage{graphicx}%
\providecommand{\U}[1]{\protect\rule{.1in}{.1in}}
\newtheorem{theorem}{Theorem}[section]
\theoremstyle{plain}

\numberwithin{equation}{section}
\begin{document}
\title[{\small rigidity of self-shrinking solutions}]{{\small A pointwise approach to rigidity of almost graphical self-shrinking
solutions of mean curvature flows}}
\author{Dongsheng LI}
\address{Dongsheng Li: School of Mathematics and Statistics\\
Xi'an Jiaotong University \\
Xi'an 710049, China}
\email{lidsh@mail.xjtu.edu.cn}
\author{Yingfeng XU}
\address{School of Mathematics and Statistics\\
Xi'an Jiaotong University \\
Xi'an 710049, China}
\email{xuyingfeng@stu.xjtu.edu.cn}
\author{Yu YUAN}
\address{Department of Mathematics, Box 354350\\
University of Washington\\
Seattle, WA 98195}
\email{yuan@math.washington.edu}

\begin{abstract}
We prove rigidity of any properly immersed noncompact Lagrangian shrinker with
single valued Lagrangian angle for Lagrangian mean curvature flows. Our
pointwise approach also provides an elementary proof to the known rigidity
results for graphical and almost graphical shrinkers of mean curvature flows.

\end{abstract}
\date{\today}
\maketitle

\section{\bigskip Introduction}

In this note, we prove the following

\begin{theorem}
If $u\left(  x\right)  $ is a smooth solution to the potential equation for
Lagrangian shrinker $\left(  x,Du\left(  x\right)  \right)  \subset
\mathbb{R}^{n}\times\mathbb{R}^{n}$
\begin{equation}
\Theta=\sum_{i}^{n}\arctan\lambda_{i}=\frac{1}{2}x\cdot Du\left(  x\right)
-u\left(  x\right)  \label{Epoten}%
\end{equation}
on (bounded or unbounded) domain $\Omega\subset\mathbb{R}^{n}$ such that
$\left\vert Du\left(  x\right)  \right\vert =\infty$ on the boundary
$\partial\Omega,$ where $\lambda_{i}s$ are the eigenvalues of $D^{2}u,$ then
the Lagrangian shrinker is a plane over $\Omega=\mathbb{R}^{n}.$ \newline More
generally, if $L^{n}$ is a smooth, properly immersed (extrinsically complete),
and noncompact Lagrangian shrinker in $\mathbb{R}^{n}\times\mathbb{R}^{n},$
where the Lagrangian angle $\Theta$ is a single valued function (zero Maslov
class), then $L^{n}$ is a Lagrangian plane.
\end{theorem}

Our pointwise approach to the shrinkers of Lagrangian mean curvature flows
also provides a short proof for the rigidity of codimension one graphical and
almost graphical shrinkers of mean curvature flows, which have been done via
integral ways by Wang [W] and Ding-Xin-Yang respectively [DXY].

\begin{theorem}
[{[W]}]Every smooth entire graphical self-shrinking hypersurface of the mean
curvature flow must be a plane.
\end{theorem}

\begin{theorem}
[{[DXY]}]Every smooth, almost graphical, properly immersed (extrinsically
complete), and noncompact self-shrinking (oriented) hypersurface of the mean
curvature flow must be a plane or a cylinder with cross section being a
self-shrinker in one lower dimensional space.
\end{theorem}

Here \textquotedblleft almost graphical\textquotedblright\ means (one choice)
of all unit normals of the hypersurface are on the closed upper hemisphere of
the whole ambient Euclid space; \textquotedblleft properness\textquotedblright%
\ or \textquotedblleft extrinsic completeness\textquotedblright\ means the
distance from every point of shrinker boundary to the origin is infinite.

Self-shrinking solutions arise naturally at a minimum-blowing-up-rate or type
I singularity from Huisken's monotonicity formula [H] for mean curvature
flows. These are immersions $F\left(  p,t\right)  :\Sigma\times\left(
-\infty,0\right)  \rightarrow\mathbb{R}^{n+k}$ which deform homothetically
$F\left(  \Sigma,t\right)  =\sqrt{-t}F\left(  \Sigma,-1\right)  $ under the
mean curvature flow equation%
\[
\left(  F_{t}\right)  ^{\perp}=\triangle_{g}F.
\]
Here $\left(  \ \ \right)  ^{\perp}$ is the normal component of the vector
$\left(  \ \ \right)  $ and $\triangle_{g}F$ equals the mean curvature of
$F\left(  \Sigma,t\right)  $ with $g$ being the induced metric. Equivalently,
the self-shrinker or shrinker $\Sigma=F\left(  p,-1\right)  $ satisfies%
\[
\triangle_{g}F=-\frac{1}{2}F^{\perp}.
\]

When shrinker $\Sigma$ is a codimension $k$ graph $\left(  x,f\left(
x\right)  \right)  \subset\mathbb{R}^{n}\times\mathbb{R}^{k},$ for the profile
$f\left(  x\right)  $ of the shrinking solution $\left(  x,\sqrt{-t}f\left(
x/\sqrt{-t}\right)  \right)  ,$ the above self-shrinking equation also takes
the following non-divergence as well as divergence form%
\begin{gather*}
g^{ij}D_{ij}f=\frac{1}{2}\left[  x\cdot Df\left(  x\right)  -f\left(
x\right)  \right]  ,\\
\triangle_{g}f^{\alpha}=\frac{1}{2}\left(  \left\langle F,\nabla_{g}f^{\alpha
}\right\rangle -f^{\alpha}\right)
\end{gather*}
for $\alpha=1,\cdots,k.$ The equivalence of the two forms comes from a simple
identity on the shrinker $\Sigma=F\left(  p,-1\right)  $%
\begin{equation}
g^{ij}D_{ij}-\frac{1}{2}x\cdot D=\triangle_{g}-\frac{1}{2}\left\langle
F,\nabla_{g}\right\rangle .\label{Ndiv=div}%
\end{equation}

When shrinker $\Sigma$ is a Lagrangian or \textquotedblleft
gradient\textquotedblright\ graph $L=\left(  x,Du\left(  x\right)  \right)
\subset\mathbb{R}^{n}\times\mathbb{R}^{n},$ the potential equation
(\ref{Epoten}) is revealed by integrating the non-divergence equation. For
each oriented tangent plane to any Lagrangian submanifold $L^{n}%
\subset\mathbb{R}^{n}\times\mathbb{R}^{n},$ there are $n$ canonical angles up
to a multiple of $2\pi$ formed with the $x$-plane $\mathbb{R}^{n},$ the sum of
those angles is called the Lagrangian angle. For example, when a Lagrangian
submanifold is a graph over $x$-space, it must be a \textquotedblleft
gradient\textquotedblright\ one $\left(  x,Du\left(  x\right)  \right)  ,$ the
Lagrangian angle is a single valued function $\Theta=\sum_{i}^{n}%
\arctan\lambda_{i}\left(  D^{2}u\right)  ;$ for circle $\left(  \cos
\theta,\sin\theta\right)  $ on $\mathbb{R}^{1}\times\mathbb{R}^{1},$ its
Lagrangian angle is multiple valued $\theta+\pi/2.$

Our current work grew out of an attempt on the rigidity issue for
extrinsically complete graphical shrinkers defined on bounded domains. We are
grateful to Tom Ilmanen for this question. Our resolution requires one to
exploit both vertical and horizontal parts of position vector of shrinkers for
a \textquotedblleft full\textquotedblright\ barrier, instead of just the
horizontal part as in the joint work [CCY] with Chau and Chen.

Heuristically, our argument goes as follows. A geometric quantity, which is
the corresponding cosine to the slope of the codimension one almost graphical
shrinker, or the Lagrangian angle of the Lagrangian shrinker, satisfies an
elliptic equation with self-similar term (Step 1s). This amplifying force term
forces the geometric quantity to go up near infinity by the \textquotedblleft
full\textquotedblright\ barrier (Step 2s). Hence the geometric quantity is
constant by the strong minimum principle, and in turn the rigidity follows
from, the second fundamental form term in the cosine equation for the
codimension one shrinkers (Step 3 of Section 2), and the quadratic excess
terms in the potential equation for the Lagrangian shrinkers (Step 3 of
Section 3).

\section{Proof Theorem 1.2 and 1.3}

Step 1. Starting from the self-shrinking equation, a simple calculation [EH,
p.471] shows that the cosine of the angle between a unit normal $N$ of the
oriented immersed shrinker $\Sigma$ and a fixed direction $e_{n+1}=\left(
0,\cdots,0,1\right)  $ in $\mathbb{R}^{n}\times\mathbb{R}^{1},$ the
nonnegative $w=\left\langle N,e_{n+1}\right\rangle $ satisfies%
\[
\mathcal{L}w=\bigtriangleup_{g}w-\frac{1}{2}\left\langle X,\nabla
_{g}w\right\rangle =-\left\vert A\right\vert ^{2}w\leq0,
\]
where $\left\vert A\right\vert $ denotes the norm of the second fundamental
form of the hypersurface $\Sigma;$ moreover, a straightforward calculation
shows that the distance $\left\vert X\right\vert $ from any point on the
shrinker $\Sigma^{n}$ to the origin satisfies%
\[
\mathcal{L}\left\vert X\right\vert ^{2}=\bigtriangleup_{g}\left\vert
X\right\vert ^{2}-\frac{1}{2}\left\langle X,\nabla_{g}\left\vert X\right\vert
^{2}\right\rangle =2n-\left\vert X\right\vert ^{2}.
\]

Step 2. Based on $\left\vert X\right\vert ^{2}$ we construct a barrier to
force $w$ to attains its global minimum at a finite point of $\Sigma.$ Set%
\[
b=-\varepsilon\left(  \left\vert X\right\vert ^{2}-K^{2}\right)
+\min_{\mathfrak{\bar{B}}_{K}\left(  0\right)  \cap\Sigma}w,
\]
where $\varepsilon$ is any fixed small positive number and $\mathfrak{\bar{B}%
}_{K}\left(  0\right)  $ is the ball in $\mathbb{R}^{n+1}$ centered at the
origin with radius $K$ such that $K\geq\sqrt{2n}$ and $\mathfrak{B}_{K}\left(
0\right)  \cap\Sigma$ is not empty. Now $\mathcal{L}b=-\varepsilon\left(
2n-\left\vert X\right\vert ^{2}\right)  \geq0$ on each unbounded component of
$\Sigma\backslash\mathfrak{\bar{B}}_{K}\left(  0\right)  ,$ and on the
infinite and finite boundary of each such component $w\geq b.$ Here we used
the properness of $\Sigma,$ or $\left\vert X\right\vert =\infty$ at the
(infinite) boundary of $\Sigma$ for the boundary comparison. By the comparison
principle $w\geq b$ on all those unbounded components of $\Sigma
\backslash\mathfrak{\bar{B}}_{K}\left(  0\right)  .$ By letting $\varepsilon$
go to zero, we then conclude that $w$ achieves its global minimum at a finite
point on $\Sigma$ (could be outside $\mathfrak{\bar{B}}_{K}\left(  0\right)
$). The strong minimum principle implies that $w$ is a constant.

\textbf{Remark.} In the case of the shrinker being an entire graph
$\Sigma=\left(  x,f\left(  x\right)  \right)  ,$ the argument in this Step 2
is \textquotedblleft cleaner\textquotedblright.

\bigskip

Step 3. If constant $w>0,$ then by the equation for $w,$ one sees that
$\left\vert A\right\vert =0$ and the almost graphical shrinker is a plane. If
constant $w\equiv0,$ then vertical vector $\left(  0,\cdots,0,1\right)  $ is
tangent to $\Sigma$ everywhere, and in turn the almost graphical shrinker is a
cylinder $\Sigma^{n}=\bar{\Sigma}^{n-1}\times\mathbb{R}^{1}$ with $\bar
{\Sigma}^{n-1}$ being a shrinker in $\mathbb{R}^{n}.$

\section{Proof of Theorem 1.1}

Step 1. When the Lagrangian shrinker is locally a graph $L=\left(  x,Du\left(
x\right)  \right)  ,$ as calculated in [CCY, p.232], we have the equation for
$\Theta$%
\[
g^{ij}D_{ij}\Theta-\frac{1}{2}x\cdot D\Theta=0.
\]
Because of (\ref{Ndiv=div}), $\Theta$ also satisfies a divergence equation%
\[
\mathcal{L}\Theta=\bigtriangleup_{g}\Theta-\frac{1}{2}\left\langle
X,\nabla_{g}\Theta\right\rangle =0
\]
with $X$ being the position vector of $L$ in $\mathbb{R}^{n}\times
\mathbb{R}^{n}.$ Note this divergence operator $\mathcal{L}$ is invariant
under any parametrization of $L.$ Again, a straightforward calculation shows
\[
\mathcal{L}\left\vert X\right\vert ^{2}=\bigtriangleup_{g}\left\vert
X\right\vert ^{2}-\frac{1}{2}\left\langle X,\nabla_{g}\left\vert X\right\vert
^{2}\right\rangle =2n-\left\vert X\right\vert ^{2}.
\]

\bigskip

Step 2. We take the same barrier%
\[
b=-\varepsilon\left(  \left\vert X\right\vert ^{2}-K^{2}\right)
+\min_{\mathfrak{\bar{B}}_{K}\left(  0\right)  \cap L}\Theta,
\]
where $\varepsilon$ is any fixed small positive number and $\mathfrak{\bar{B}%
}_{K}\left(  0\right)  $ is the ball in $\mathbb{R}^{n}\times\mathbb{R}^{n}$
centered at the origin with radius $K$ such that $K\geq\sqrt{2n}$ and
$\mathfrak{B}_{K}\left(  0\right)  \cap L$ is not empty. Now $\mathcal{L}%
b=-\varepsilon\left(  2n-\left\vert X\right\vert ^{2}\right)  \geq0$ on each
unbounded component of $L\backslash\mathfrak{\bar{B}}_{K}\left(  0\right)  ,$
and on the infinite and finite boundary of each such component $\mathsf{\Theta}\geq b.$ Here
we used the properness of $L,$ or $\left\vert X\right\vert =\infty$ at the
(infinite) boundary of $L$ for the boundary comparison. In the case $L=\left(
x,Du\left(  x\right)  \right)  \subset\Omega\times\mathbb{R}^{n},$ $\left\vert
X\right\vert =\infty$ is because of either $\left\vert Du\left(  x\right)
\right\vert =\infty$ or $\left\vert x\right\vert =\infty$ on the (infinite)
boundary of $L.$ By the comparison principle $\mathsf{\Theta}\geq b$ on all
those unbounded components of $L\backslash\mathfrak{\bar{B}}_{K}\left(
0\right)  .$ By letting $\varepsilon$ go to zero, we then conclude that
$\mathsf{\Theta}$ achieves its global minimum at a finite point on $L$ (could
be outside $\mathfrak{\bar{B}}_{K}\left(  0\right)  $). The strong minimum
principle implies that $\mathsf{\Theta}$ is a constant.

\bigskip

Step 3. We go to the potential equation (\ref{Epoten}) to capture the flatness
of the Lagrangian shrinker. Near a closest point $P$ on $L$ to the origin, one
can represent $L$ as a \textquotedblleft gradient\textquotedblright\ graph
$\left(  x,Dv\left(  x\right)  \right)  $ over the Lagrangian plane through
the origin and parallel to the tangent plane of $L$ at $P.$ Here we use the
abused notation $x$ for the coordinates on the ground, or the Lagrangian
plane. Because of the constancy of the Lagrangian angle $\Theta,$ the
potential equation for $v$ becomes%
\[
c=\frac{1}{2}x\cdot Dv\left(  x\right)  -v\left(  x\right)  .
\]
This $c$ may differ from the constant $\Theta$ by another one due to a
possible $U\left(  n\right)  $ coordinate rotation in $\mathbb{C}%
^{n}=\mathbb{R}^{n}\times\mathbb{R}^{n}.$ Euler's homogeneous function theorem
implies that the smooth function $v\left(  x\right)  +c$ around $x=0$ is a
polynomial of degree two. We immediately see that near $P,$ $L=\left(
x,Dv\left(  x\right)  \right)  $ is a piece of the above \textquotedblleft
ground\textquotedblright\ plane. Because of the analyticity of \ $L,$ as the
potential equation is analytic, we conclude that Lagrangian shrinker $L$ is a
Lagrangian plane, and over $\Omega=\mathbb{R}^{n}$ in the graphical case.

\bigskip

\textbf{Acknowledgments. }The first and the second authors are partially
supported by NSFC. 11671316. The third author is partially supported by an NSF grant.

\end{document}